\documentclass[12pt]{amsart}
\usepackage{graphicx}
\usepackage{amsmath}
\usepackage{amssymb}
\usepackage{a4wide}

\newcommand{\D}{\mathbb{D}}
\numberwithin{equation}{section}

\title[Estimates for integrals of derivatives of rational functions]{Estimates for integrals of derivatives \\ of rational functions in multiply \\ connected domains in the plane}
\author{A.\,D.~Baranov}
\address{Saint Petersburg State University}
\email{anton.d.baranov@gmail.com} 
\author{I.\,R.~Kayumov}
\address{Kazan Federal University}
\email{ikayumov@kpfu.ru}


\thanks{A.\,D. Baranov was supported by Ministry of Science and Higher Education 
of the Russian Federation, agreement No 075-15-2021-602.  
I.\,R. Kayumov was supported by Ministry of Science and Higher Education 
of the Russian Federation, agreement No 075-15-2019-1619.}

\begin{document}
\sloppy

\begin{abstract}
We obtain estimates for integrals of derivatives of rational functions  in multiply connected domains in the plane.
A sharp order of the growth is found for the integral of the modulus of the derivative of a finite Blaschke product in the unit disk. 
We also extend the results of E.\,P. Dolzhenko about the integrals of the derivatives of rational functions to a wider class of domains, namely, to domains bounded by rectifiable curves without zero interior angles, and show the sharpness of the obtained results. 
\end{abstract}

\maketitle

\section{Introduction}

70 years ago S.\,N.~Mergelyan \cite{Mer} showed that there exists a bounded analytic function $f$ in the disk
$\D=\{|z|<1\}$ such that 
$$
I(f):=\int_\D |f'(z)|\, dA(z) = \infty,
$$
where $dA(z) = \frac{1}{\pi} dxdy$, $z=x+iy$.

This problem was further investigated by W.~Rudin \cite{Rud} who constructed an infinite Blaschke product
$$
B(z)=\prod_{k=1}^\infty \frac{|z_k|}{z_k}\frac{z_k-z}{1-\overline{z_k} z}
$$
such that
$I(B)=\infty$  and, moreover, $\int_{0}^1 |B'(re^{i\theta})| dr = \infty$ for a.e. $\theta\in[0, 2\pi]$.
A similar, but more explicit example, was given by G. Piranian \cite{Pir}.

It is then natural to ask what happens if $B$ is a finite Blaschke product of degree $n$?
It is obvious, that, for any fixed $n$, the quantity $I(B)$ is bounded, but it cannot be uniformly bounded with respect to 
$n$, since any bounded function in $\D$ is a locally uniform limit of finite Blaschke products. We find the sharp order of growth
for such integrals. Namely, we have
\medskip
\\
{\bf Theorem  1.} {\it
	Let $B$ be a finite Blaschke product of degree $n$. Then
	\begin{equation} \label{upper-estimate}
		I(B)  \leq \pi(1+\sqrt{\log n}).
	\end{equation}
On the other hand, there exists an absolute constant $c>0 $ such that for any $n \in \mathbb{N}$
	there exists a finite Blaschke product of degree $n$ satisfying $I(B) \ge c(1+\sqrt{\log n})$.}
\medskip

The proof of sharpness of this inequality is based on subtle results of N.\,G.~Makarov \cite{Mak} and 
R.~Ba\~nuelos and C.\,N.~Moore \cite{BanMoo} on boundary behaviour of functions from the Bloch space. 

It should be noted that there exists a vast literature dealing with the membership
of the derivatives of the Blaschke products to various functional spaces, e.g., Bergman-type spaces
(see \cite{av, prot1, prot2} and the references therein). However, most of these results
concern infinite products and the conditions are formulated in terms of their zeros. 

Since a Blaschke product is a bounded rational function in the unit disk, the problem about the
estimates of the derivatives of Blaschke products is related to a more general question 
about the integrals  of bounded rational functions studied for the first time by E.\,P.~Dolzhenko \cite{Dol}
for sufficiently smooth domains. We will say that a curve belongs to the class K if it is a closed
Jordan curve with continuous curvature $k(s)$ satisfying a H\"older condition as the function of the arc length $s$. 
Let $G$ be a finitely connected domain whose boundary curves belong to the class K.
Assume that $1\le p \le 2$ and let $R$ be a rational function of degree at most $n$ with the poles 
outside $\overline{G}$. Dolzhenko \cite[Theorem 2.2]{Dol} showed that there exists 
a constant $C$ depending only on the domain $G$ and on $p$ such that
\begin{equation}
	\label{d1} 
	\int_{G}|R'(w)|^p\, dA(w) \leq C n^{p-1} \|R\|_{H^\infty(G)}^p, \qquad p \in (1,2],
\end{equation}
\begin{equation}
	\label{d2} 
	\int_{G}|R'(w)|\, dA(w) \leq C \ln (n+1) \|R\|_{H^\infty(G)}.
\end{equation}
Here we denote by  $H^\infty(G)$ the space of all bounded analytic functions in $G$, 
and $\|f\|_{H^\infty(G)} = \sup_{w\in G} |f(w)|$.

Later, inequalities for the derivatives of rational functions (mainly in the disk) were studied in the papers by 
V.\,V.~Peller \cite{pel}, S.~Semmes \cite{sem}, A.\,A.~Pekarskii \cite{pek, pek1}, V.\,I.~Danchenko 
\cite{dan0, dan} and by many other authors (see, e.g., \cite{dyn1, dyn2, bz0, bz1, bz2}). 
A short proof of the Dolzhenko inequalities for the case of the disk when the $H^\infty$-norm is replaced by the weaker
$BMOA$-norm  can be found in \cite{bz1}.

In the present article the inequalities \eqref{d1} and \eqref{d2} are proved under substantially weaker 
restrictions on the domain, namely, under the condition that the domain has no zero interior angles 
(more precisely, for the John class domains -- see the definition in  \S 3).
\medskip
\\
{\bf Theorem 2.} {\it Let $G$ be a finitely connected John domain with the rectifiable boundary 
and let $1\le p \le 2$. Then there exists a constant $C>0$, depending on the domain $G$ and on $p$, such that
for any rational function $R$ of degree at most $n$ the inequalities \eqref{d1} and \eqref{d2} hold.} 
\medskip

The sharpness of \eqref{d1} is seen already on the simplest example of the function 
$R(z) = z^n$ in the disk (obviously, we can consider polynomials as a special case of rational
functions with the pole at infinity). The question about sharpness of the estimate \eqref{d2} 
in the conditions of Theorem 2  remains open. However, it turns out that under some additional regularity 
of the domain $G$ inequality \eqref{d2} can be improved. 
\medskip
\\
{\bf Theorem 3.}
{\it Let $G$ be a simply connected domain such that 
$\varphi' \in H^2$, where $\varphi$ is the conformal map of the disk $\mathbb{D}$ onto $G$.
	Then there exists a constant $C>0$ depending on the domain $G$ such that  
	for any rational function $R$ of degree at most $n$ one has
	\begin{equation} 
		\label{g1}
		\int_{G}|R'(w)|\, dA(w) \leq C \sqrt{\ln (n+1)} \|R\|_{H^\infty(G)}.
\end{equation}}

As follows from Theorem 1, the dependence on $n$ in this inequality is sharp.  
\medskip

Finally, we give a statement for the case $p> 2$. Here $q$ is the conjugate exponent, i.e., $1/p+1/q=1$.
\medskip
\\ 
{\bf Theorem 4. } { \it Let  $G$ be a bounded simply connected domain and put 
	$G_\rho = \{ z\in G: {\rm dist}\,(z, \partial G) >\rho\}$. Then for any rational function $R$ of degree at most $n$
   and $p>2$ one has 
	\begin{equation} 
		\label{g2}
		\|R'\|_{A^p(G_\rho)}  = \bigg(\int_{G_\rho}|R'(w)|^p\, dA(w)\bigg)^{1/p} \le  n^{1/p} \rho^{1/p-1/q} \|R\|_{H^\infty(G)}.
\end{equation} }

In \cite{Dol} inequality \eqref{g2} was established for domains of the class K, but,
as our result shows, no restrictions on the regularity of the domain are required. 

In \cite{bz2} another generalization of the Dolzhenko inequality for the case $p>2$ was obtained 
for the rational functions in the disk. Let $R$ be a rational function  of degree at most $n$
whose poles lie in the complement of the disk $\{|z| <1+\rho\}$. The following inequality 
follows directly from \cite[Theorem 8.2]{bz1}:  
$$
\|R'\|_{A^p(\D)} \le  C(p) n^{1/q} \rho^{1/p-1/q} \|R\|_{BMOA};
$$
here $BMOA$ denotes the analytic space of functions of bounded mean oscillation in the disk.
It is interesting to note that the dependence on $n$ in Theorem 4 is substantially weaker
(since in this case the function $R$ is assumed to be bounded on a larger set). 

In \S 5 more general inequalities are obtained for the weighted norms of the derivatives of rational functions,
where the weight is given as some power of the distance to the boundary of the domain. 

A suitable toolbox for the study of such inequalities is provided by the theory of the Hardy spaces. For $p>0$
the Hardy  space $H^p$ is the set of all analytic functions in $\mathbb{D}$ satisfying 
$\|f\|_{H^p}<\infty$, where
$$
\|f\|^p_{H^p}:=\sup_{0<r<1}\frac{1}{2\pi}\int_0^{2\pi}|f(re^{it})|^p dt.
$$
Note that for $p \ge 1$ the last quantity defines a norm with respect to which $H^p$ is a Banach space.  
\bigskip


\section{Estimate for the integral of the modulus of the derivative of a finite Blaschke product }

In the proof of Theorem 1 we will use the following simple lemma.
\medskip
\\
{\bf Lemma 1.} {\it Assume that the function $g(z) = \sum_{k=0}^\infty b_k z^k$ 
    is analytic in the disk $\mathbb{D}$. If $\|g\|_\infty \le 1$
	and $p(z) = \sum_{k=0}^n b_k z^k$, $n\ge 2$, then there exists an absolute constant $C_0$ such that
	$$
	|p(z)| \le C_0, \qquad |z|\le 1-2\log n/ n,
	$$
	and
	$$
	|g'(z) - p'(z)| \le C_0, \qquad |z|\le 1-2\log n/ n.
	$$}

\noindent
{\bf Proof.} Since $|b_k| \le 1$ and $|z|^k \le 1/n^2$ for $|z|\le 1-2\log n/ n$ and
$k\ge n$, the function $|g(z) - p(z)|$ admits a uniform estimate for $|z|\le 1-2\log n/ n$, 
and the first estimate follows. 

Clearly, 
$$
\sum_{k=n}^\infty (k+1)|b_{k+1} z^k| \le \frac{|z|^n}{(1-|z|)^2} +\frac{n|z|^n}{1-|z|} \le const, \qquad |z|\le 1-2\log n/ n,
$$
and the second inequality is proved. 
\bigskip

\noindent
{\bf Proof of Theorem.}  We use the following well-known facts:
$$
\int_0^{2\pi} |B'(re^{it})| dt  \le 2 \pi n, \qquad r \in[0,1],
$$
for any finite Blaschke product of degree at most $n$ and
\begin{equation}
\label{area}
\int_\D |f'(z)|^2 (1-|z|^2)\, dA(z)  = 
\sum_{n=1}^\infty \frac{n}{n+1}|a_n|^2 \le \|f\|^2_{H^2}
\end{equation}
for any function $f(z) = \sum_{n\ge 0} a_n z^n$ in the Hardy space $H^2$.

Let $s \in [0,1]$. We have
$$
\int_{\{s<|z| <1\}} |B'(z)|\, dA(z)= \int_0^{2\pi} \int_s^1 |B'(re^{it})| r drdt \leq 2\pi n \int_s^1  r dr=\pi n (1-s^2).
$$
In the remaining part we apply the Cauchy--Schwarz inequality:
$$
	\int_{\{0 <|z| \le s\}} |B'(z)|\, dA(z) 
$$
$$
	\le \bigg( \int_{\{0 <|z| \le s\}} 
   |B'(z)|^2 (1-|z|^2)\, dA(z) \bigg)^{1/2} \bigg( \int_{\{0 <|z| \le s\}} \frac{dA(z)}{1-|z|^2}
	\bigg)^{1/2}
$$
$$
	 \leq	\sqrt{\pi}\sqrt{2\pi \int_0^s \frac{rdr}{1-r^2}} = \pi \sqrt{ \log\frac{1}{1-s^2}}.
$$
Thus,
$$
I(B) \leq \pi n (1-s^2)+\pi \sqrt{ \log\frac{1}{1-s^2}}.
$$ 
Taking $s^2=1-1/n$ we obtain \eqref{upper-estimate}.

The estimate from below can be obtained by the methods based on the 
Makarov law of the iterated logarithm \cite{Mak}. Recall that the Bloch class $\mathcal{B}$ consists of functions 
analytic in $\D$ with finite seminorm $ \|f\|_{\mathcal{B}}  = \sup_{z\in \D} (1-|z|^2) |f'(z)|$.
In \cite{BanMoo} R.~Ba\~nuelos and C.\,N.~Moore, answering a question of N.\,G.~Makarov and F.~Przytycki,
constructed a function $f(z) = \sum_{k=1}^\infty a_k z^k$ in the Bloch class such that its asymptotic entropy  
admits the lower bound
$$
\liminf_{r\to 1-} \frac{\sum_{k=1}^\infty |a_k|^2 r^{2k}}{\log\frac{1}{1-r}} >0,
$$
whereas for all $\zeta$ with $|\zeta| =1$ one has
$$
\limsup_{r\to 1-} \frac{f(r \zeta)}{\sqrt{\log\frac{1}{1-r} \log\log\log \frac{1}{1-r}}} =0.
$$
Moreover, in \cite[p. 852--853]{BanMoo} a sequence of polynomials  
$$
p_n(z) = \sum_{k=4}^{4^{n+1} -1} a_k z^k = \sum_{j=1}^{n} b_j(z), \qquad \text{where} \qquad 
b_j(z) = \sum_{k=4^j}^{4^{j+1} -1} a_k z^k,
$$ 
is constructed such that $\|b_j\|_\infty \le 1$,
$$
\sum_{k=1}^{4^{n+1} -1} |a_k|^2 \ge c \log m, 
\qquad  
\|p_n\|_{H^\infty} \le  C\sqrt{\log m},
$$
where  $m = {\rm deg}\, p_n  = 4^{n+1} -1$ and $C, c>0$ are some absolute positive constants.  

It is not difficult to deduce from $\|b_j\|_\infty \le 1$ that $\sup_n \|p_n\|_{\mathcal{B}} <\infty$. 
Indeed, by the Schwarz lemma $|b_j(rz)| \le r^{4^j}$, whereas, by the classical Bernstein inequality, 
$|b_j' (rz)| \le 4^j r^{4^j}$. It is well known (and easy to show) that 
$\sum_{j\ge 1} 4^j r^{4^j} \le C_1/(1-r^2)$ for some constant $C_1>0$ and, thus,  
$\sup_n \|p_n\|_{\mathcal{B}} <\infty$. Without loss of generality we may assume that $\|p_n\|_{\mathcal{B}} \le 1$.

Let $r=1-1/m$. Then, for some absolute constants $C', c'>0$,
$$
c' \log\frac{1}{1-r} \leq \sum_{k=1}^m |a_k|^2 r^{2k} \leq 2 \int_{|z|<r}|p_n'(z)|^2(1-|z|^2)\,dA(z)
\leq C'\int_{|z|<r}|p_n'(z)|\,dA(z).
$$
Now put $q_n = p_n/(C\sqrt{\log m})$. Then $\|q_n\|_\infty \le 1$ and
$$
\int_{|z|<1-1/m}|q_n'(z)|\, dA(z)  \ge c_1\sqrt{\log m}
$$
for some absolute constant $c_1>0$. Since
$$
\int_{1-2\log m/m < |z|<1-1/m}|q_n'(z)| \, dA(z) \leq \int_{1-2\log m /m < |z|<1-1/m} \frac{dA(z)}{1-|z|^2}=O(\log\log m),
$$
we have, for some $c_2>0$,
\begin{equation}
	\label{eq1}
	\int_{|z|<1-2\log m/m}|q_n'(z)|\,dA(z)  \ge c_2 \sqrt{\log m}.
\end{equation}

Let $B$ be a Blaschke product of degree at most $m+1$ such that its first $m$ Taylor coefficients coincide with respective coefficients of the polynomial $q_n$. By Lemma 1, 
$$
|q_n'(z)-B'(z)| \leq 2C_0, \qquad |z| \leq 1-2\log m/m.
$$
Hence, it follows from \eqref{eq1} that
$$
\int_{\D}|B'(z)|\, dA(z) \ge c_3\sqrt{\log m}
$$
for some absolute constant $c_3$.  Theorem 1 is proved.
\bigskip

\section{Estimates for integrals of rational functions}

Recall that a finitely connected domain $\Omega$ is said to be a {\it John domain} 
if there exists a constant $C>0$  such that any points $a,b\in \Omega$ 
can be connected by a curve $\gamma$ in $\Omega$ with the following property:
for any $x\in\gamma$,
$$
\min \big( {\rm diam}\,\gamma(a,x), {\rm diam}\,\gamma(x,b) \big) \le C {\rm dist}\, (x, \partial \Omega).
$$
Here $\gamma (a,x)$ and $\gamma (x,b)$ denote the corresponding subarcs of $\gamma$. 
For equivalent definitions and properties of John domains see, e.g., \cite{MaSa, Pom}. 
Essentially, this definition means that the domain has no zero interior angles. In particular,
a domain is a John domain if it satisfies the cone conditon: one can touch each boundary point
from inside of the domain by some sufficiently small triangle with fixed angles. 

In what follows we will essentially use the following property of simply connected John domains: 
if $\varphi$ is the conformal map of $\D$ onto a simply connected John domain, then
\begin{equation}
	\label{hold}
	|\varphi'(z)| \le \frac{C}{(1-|z|)^{\alpha}}
\end{equation}
for some $\alpha\in(0,1)$ and  $C>0$ (see \cite[pp. 96--100]{Pom}).

In the proof of Theorem 2 we will use the following simple lemma: 
\medskip
\\
{\bf Lemma 2.} {\it Let $g$ be a bounded and at most $n$-valent function in $\D$. Then
	$$
	\frac{1}{2\pi} \int_0^{2\pi} |g'(re^{it})|^2 dt \le \frac{n}{1-r} \|g\|^2_{H^\infty(\D)}.
	$$}

\noindent
{\bf Proof.} Let $g(z) = \sum_{k=0}^\infty a_kz^k$. Then
$$
\frac{1}{2\pi} \int_0^{2\pi} |g'(re^{it})|^2 dt = \sum_{k=1}^\infty k^2|a_k|^2 r^{2k} \leq
\frac{1}{1-r} \sum_{k=1}^\infty k|a_k|^2.
$$
Here we used the elementary inequality $kr^{k-1}(1-r) \leq 1$, $k\ge 1$, $r \in [0,1)$.
Since $g$ is at most $n$-valent in $\mathbb{D}$, we have, in virtue of the classical area theorem,
$$
\sum_{k=1}^\infty k|a_k|^2 \leq  n \|g\|^2_{H^\infty(\D)}.
$$

\noindent
{\bf Proof of Theorem 2.}
Assume first that the domain $G$ is bounded. Without loss of generality we may assume that $G$ 
is a simply connected domain with a rectifiable boundary. Indeed, making smooth cuts (and controlling the angles) one can easily 
represent our domain as a finite union of simply connected John domains. 

Let $w=\varphi(z)$ be a conformal map of $\D$ onto $G$.
Since the boundary of $G$ is rectifiable, we have $\varphi' \in H^1$. Also, $\varphi'$ satisfies inequality \eqref{hold}.

By the change of the variable, 
$$
	\int_{G}|R'(w)|^p\, dA(w) = \int_{\mathbb{D}}|R'(\varphi(z))|^p|\varphi'(z)|^2\, dA(z)= 
	\int_{\mathbb{D}} |(R\circ \varphi)'(z)|^p  |\varphi'(z)|^{2-p} dA(z).
$$ 
It is obvious, that for $p=2$ the last integral does not exceed $n \|R\|_{H^\infty(G)}^2$, since the function $R\circ\varphi$ 
is at most $n$-valent  in the disk $\mathbb{D}$.

Let us split the last integral into the integrals over the set $\{|z| \le r_n\}$ and over the set $\{r_n < |z| < 1\}$, 
where $r_n = 1- \frac{1}{(n+1)^{K}}$ and the number  $K>0$ is to be chosen later. 
\medskip
\\
{\bf Estimate of the integral over the set $\{|z| \le r_n\}$.} Let $M = \|R\|_{H^\infty(G)}$. 
Set
$$
J := \int_{\{|z|\le r_n\}} |(R\circ \varphi)'(z)|^p  |\varphi'(z)|^{2-p} dA(z). 
$$

For $p=1$ we use the estimate $(1-|z|^2) |(R\circ \varphi)'(z)| \le M$. We have
$$
\begin{aligned}
	J  \le \frac{M}{\pi} \int_0^{r_n} & \frac{1}{1-r}\int_0^{2\pi} |\varphi'(re^{it})| dt dr 
	\\ &
	\le  2 \|\varphi'\|_{H^1} M 
	\int_0^{r_n} \frac{dr}{1-r} = 2K \|\varphi'\|_{H^1} \log(n+1) M. 
\end{aligned}
$$

In the case $1<p<2$ consider separately the integrals over the sets  
$\{|z| \le 1-\frac{1}{n+1}\}$ and ${\{1-\frac{1}{n+1} <|z| \le r_n\}}$. Since $\varphi' \in H^1$, we have
$\varphi'\in H^{2-p}$ and $\|\varphi'\|_{H^{2-p}} \le \|\varphi'\|_{H^1}$. Hence, 
$$
\begin{aligned}
	\int_{\big\{|z| \le 1- \frac{1}{n+1}\big\} }  |(R\circ \varphi)'(z)|^p  |\varphi'(z)|^{2-p} dA(z) 
	& \le 2\|\varphi'\|_{H^1}^{2-p} \int_0^{1-\frac{1}{n+1}} \frac{M^p}{(1-r)^p} dr \\
	& = 
	2\|\varphi'\|_{H^1}^{2-p} (p-1)^{-1} (n+1)^{p-1} M^p. 
\end{aligned}
$$
To estimate the integral over the set $\{1-\frac{1}{n+1} <|z| \le r_n\}$ we use the H\"older inequality 
with exponents $(2-p)^{-1}$ and $(p-1)^{-1}$:
$$
J  \le 
2 \int_{1-\frac{1}{n+1}}^{r_n} \bigg(\frac{1}{2\pi} \int_0^{2\pi} |(R\circ \varphi)'(re^{it})|^{\frac{p}{p-1}} dt\bigg)^{p-1} 
\bigg(\frac{1}{2\pi} \int_0^{2\pi} |\varphi'(re^{it})| dt\bigg)^{2-p} dr. 
$$
Using successively the inequality $(1-|z|^2) |(R\circ \varphi)'(z)| \le M$ and Lemma 2 we get 
$$
\begin{aligned}
	J &\le 2 \|\varphi'\|_{H^1}^{2-p} M^{p-2(p-1)} \int_{1-\frac{1}{n+1}}^{r_n} \frac{1}{(1-r)^{p-2(p-1)}} 
	\bigg(\frac{1}{2\pi} \int_0^{2\pi} |(R\circ \varphi)'(re^{it})|^2 dt\bigg)^{p-1}  dr  \\
	& \le  
	2 \|\varphi'\|_{H^1}^{2-p} M^p n^{p-1}
	\int_{1-\frac{1}{n+1}}^{r_n} \frac{dr}{(1-r)^{2-p +p-1}} \\
	& = 2 \|\varphi'\|_{H^1}^{2-p} M^p n^{p-1}
	\int_{1-\frac{1}{n+1}}^{1-\frac{1}{(n+1)^K}} \frac{dr}{1-r}   = 2 \|\varphi'\|_{H^1}^{2-p} K n^{p-1} M^p.
\end{aligned}
$$
\smallskip
\\
{\bf Estimate of the integral over the set $\{r_n<|z| <1\}$.} 
In this case the argument applies to all $p\in [1, 2)$.
Choose $\delta$ such that $0<\delta<1-p/2$. Then $2-p-\delta \in (0,1)$. Let
$\beta$ be the exponent conjugate to $(2-p-\delta)^{-1}$. One has, by \eqref{hold},
$$
	I  := \int_{r_n<|z|<1} |(R\circ \varphi)'(z)|^p  |\varphi'(z)|^{2-p} dA(z)
	 $$
	 $$
	\le C^\delta
	\int_{r_n<|z|<1} \frac{|(R\circ \varphi)'(z)|^p |\varphi'(z)|^{2-p-\delta}}{(1-|z|)^{\alpha\delta}} dA(z) $$$$
	 \le
	2C^\delta \int_{r_n}^1 \frac{1}{(1-r)^{\alpha\delta}} \bigg(\frac{1}{2\pi}
	\int_0^{2\pi} |(R\circ \varphi)'(re^{it})|^{p\beta} dt\bigg)^{1/\beta} 
	\bigg(\frac{1}{2\pi} \int_0^{2\pi} |\varphi'(re^{it})| dt\bigg)^{2-p-\delta} dr. 
$$
Note that it follows from $\delta<1-p/2$ that $p\beta>2$. 
Applying the estimate $(1-|z|^2) |(R\circ \varphi)'(z)| \le M$ and Lemma 2, we get
$$
\begin{aligned}
	I  & \le 
	2 C^\delta\|\varphi'\|_{H^1}^{2-p-\delta} M^{\frac{p\beta - 2}{\beta}}
	\int_{r_n}^1 \frac{1}{(1-r)^{\alpha\delta + \frac{p\beta - 2}{\beta}}} \bigg(\frac{1}{2\pi}
	\int_0^{2\pi} |(R\circ \varphi)'(re^{it})|^2 dt\bigg)^{1/\beta} dr \\
	& \le 
	2C^\delta\|\varphi'\|_{H^1}^{2-p-\delta} M^p n^{1/\beta} \int_{r_n}^1 \frac{dr}{(1-r)^{\alpha\delta + \frac{p\beta - 2}{\beta} +\frac{1}{\beta}}} \\
	& = 
	2 C^\delta\|\varphi'\|_{H^1}^{2-p-\delta} M^p n^{1/\beta} \int_{r_n}^1 \frac{dr}{(1-r)^{\alpha\delta + p -\frac{1}{\beta}}}.
\end{aligned}
$$
It remains to notice that $\alpha\delta +p- \frac{1}{\beta} = \alpha\delta +p- (1- (2-p-\delta)) = 1-(1-\alpha)\delta$, whence
\begin{equation}
	\label{pol}
	I \le 2(1-\alpha)^{-1}\delta^{-1}
   C^\delta\|\varphi'\|_{H^1}^{2-p-\delta} M^p n^{1/\beta} (1 - r_n)^{(1-\alpha)\delta}.
\end{equation}
If we fix $\delta\in (0, 1-p/2)$ and choose a sufficiently large $K$ in $r_n = 1- \frac{1}{(n+1)^K}$, 
we conclude that 
$I \le C^\delta \|\varphi'\|_{H^1}^{2-p-\delta} M^p$ (and even $o(1)$ as $n\to\infty$).
\medskip

We now consider the case when $\infty \in G$. 
It is clear that such domain can be represented as a union (possibly with an intersection) 
of the complement of a disk of sufficiently large radius with a simply connected bounded John domain. The statement is already proved for 
bounded simply connected domains, while for the complement of a disk it follows from the results of Dolzhenko cited above. 
Theorem 2 is proved.
\bigskip


\section{Proofs of Theorems 3 and 4}
\noindent
{\bf Proof of Theorem 3.}
As in the proof of Theorem we set  
$r_n = 1- \frac{1}{(n+1)^K}$, where $K>0$. Since $\varphi' \in H^2 \subset H^1$ and the condition  \eqref{hold} 
is satisfied with $\alpha = 1/2$, one can use the estimate \eqref{pol} for the integral over the set $\{r_n < |z| < 1\}$ 
established in the proof of Theorem 2. For a sufficiently large $K$ this integral is uniformly bounded over $n$ 
(and even tends to zero as $n\to\infty$). 

Thus, it suffices to estimate 
$$
J:=\int_{0<|z| \le r_n} |(R\circ \varphi)'(z)|  |\varphi'(z)| dA(z).
$$
By the Cauchy--Schwarz inequality, 
$$
\begin{aligned}
	J & \le \bigg(\int_{0<|z|\le r_n} (1-|z|) |(R\circ \varphi)'(z)|^2  dA(z) \bigg)^{1/2} 
	\bigg(\int_{0<|z| \le r_n}   \frac{|\varphi'(z)|^2}{1-|z|}dA(z) \bigg)^{1/2} \\
	& \le M \bigg( \int_{0<|z| \le r_n}   \frac{|\varphi'(z)|^2}{1-|z|} dA(z) \bigg)^{1/2} \leq
	\sqrt{2} M \|\varphi'\|_{H^2}  \bigg( \int_0^{r_n} \frac{dr}{1-r} \bigg)^{1/2} \\
	& =\sqrt{2} M \|\varphi'\|_{H^2} \sqrt{K\ln (n+1)}.
\end{aligned}
$$
Theorem 3 is proved.
\bigskip
\\
{\bf Proof of Theorem 4.}
Let  $\varphi$ be a conformal map of $\D$ onto $G$ and let $D_\rho = \varphi^{-1}(G_\rho)$.
Since $\rho \le (1-|z|^2)|\varphi'(z)|$, $z\in D_\rho$, we have
$$
\begin{aligned}
	\int_{G_\rho} |R'(\zeta)|^p dA(\zeta) & = \int_{D_\rho} |(R\circ \varphi)'(z)|^p |\varphi'(z)|^{2-p} dA(z) \\
	& \le \rho^{2-p} \int_{D_\rho} |(R\circ \varphi)'(z)|^p (1-|z|^2)^{p-2} dA(z) \\
	& \le \rho^{2-p} M^{p-2} \int_{D_\rho} |(R\circ \varphi)'(z)|^2 dA(z) \le \rho^{2-p} n M^p. 
\end{aligned}
$$
At the last step we used the fact that $R\circ \varphi$ covers each point of the disk of radius $M$ 
with multiplicity at most $n$.
Theorem 4 is proved.
\bigskip


\section{Weighted inequalities of Dolzhenko and Peller type}

As a natural generalization of the Dolzhenko inequalities one can consider  weighted integrals
of derivatives of rational functions. Similar inequalities were studied extensively in the setting of the Bergman (or Besov) 
spaces. E.g., a well-known inequality by  V.V.~Peller \cite{pel} states that for a rational function $R$ of degree $n$ 
with poles outside  $\overline{\mathbb{D}}$ one has
$$
\|R\|_{B_p^{1/p}} \le C  n^{1/p} \|R\|_{BMOA},
$$               
where $B_p^{1/p}$ is the Besov space, $p>0$, $C=C(p)$. In particular, for $1<p<\infty$,
$$
\int_\D |R'(z)|^p (1-|z|)^{p-2} dA(z) \le C n \|R\|^p_{H^\infty}.
$$
Various proofs and generalizations of this inequality can be found in \cite{sem, pek, pek1, bz1}.

Using the methods of \S 3 one can obtain more general weighted estimates where the weight equals 
to some power of the distance to the boundary. To formulate the corresponding result, we set,
for a bounded domain $G\subset \mathbb{C}$ and $z\in G$,
$$
d_G(z) : = {\rm dist}\, (z, \partial G).
$$
For $p\ge 1$, $\beta \in \mathbb{R}$ and a function $f$ analytic in $G$ put
$$
I_{p,\beta} (f) := \int_G |f'(\zeta)|^p \, d_G^\beta(\zeta)\, dA(\zeta)
$$
(in general, the quantity $I_{p,\beta} (f)$ can be infinite). We are interested in the estimates of the form   
$$
I_{p,\beta}(R) \le C\Psi(n) \|R\|^p_{H^\infty(G)}, 
$$
which hold for all  rational functions $R$ of degree at most $n$ with poles outside $\overline{G}$ and with a constant
$C$, depending on $G$, $p$ and $\beta$, but not on $n$ and $R$. Here $\Psi$ is some function depending on $n$ only. 
It is easy to see that such estimates are possible only for $\beta \ge p-2$;
it is seen already from the example $G=\D$ and rational fractions $R(\zeta) = \frac{1}{\zeta-\lambda}$ that for $\beta< p-2$ 
the integral $I_{p,\beta} (f)$ does not admit the estimate depending only on $n$, one has also to take into account 
the distance from the poles of $R$ to  $\partial G$. 

To simplify the notations, in what follows we write $X(R,n) \lesssim Y(R,n)$, if $X(R, n)\le CY(R,n)$
with a constant $C$,  depending only from $G$, $p$ and   $\beta$, but not on $n$ and $R$.
\medskip

\noindent
{\bf Theorem 5.} {\it Let $G$ be a simply connected bounded domain, $\varphi$ is the conformal map
of $\mathbb{D}$ onto $G$, 	$p\ge 1$, $\beta \ge p-2$.  The following estimates hold true.
	\smallskip
	
	1. If $\beta > p-1$ and $\varphi'\in H^{\gamma}$ for some $\gamma >1$, then 
   $I_{p,\beta}(R) \lesssim \|R\|^p_{H^\infty(G)}$, i.e., the dependence on $n$ disappears. 
	\smallskip
	
	2. If $\beta  = p-1$, $1\le p <2$ and $\varphi'\in H^{\frac{2}{2-p}}$, then
	$$
	I_{p,\beta}(R) \lesssim \big( \log n\big)^{1- \frac{p}{2}} \|R\|^p_{H^\infty(G)}.
	$$
	If $\beta  = p-1$, $p \ge 2$ and $\varphi'\in H^\infty$, then $I_{p,\beta}(R) \lesssim \|R\|^p_{H^\infty(G)}$.
	\smallskip
	
   	3. If $p-2 \le \beta <p-1$, $p\ge 2$, $\varphi' \in H^1$ and $G$ is a John domain, then
	$$
	I_{p,\beta}(R) \lesssim n^{p-1-\beta} \|R\|^p_{H^\infty(G)}.
	$$   }
\medskip

The dependence on $n$ in the inequalities in Theorem 5 is sharp already for the case of the unit disk. 
In statement 3 the optimal growth is attained on $R(z) = z^n$, while the sharpness of the inequality
in statement 2 can be shown by considering the  Ba\~nuelos--Moore construction (as $R$ one can take a polynomial 
or a Blaschke product). 

Note that statement 3 of the theorem does not cover the case $p-2 \le \beta <p-1$ and $1<p<2$. 
In this case it would be sufficient to prove the following analogue of Peller's inequality: 
$$
  \int_\D |(R\circ\varphi)'(z)|^p (1-|z|)^{p-2} dA(z) \lesssim  n\|R\|_{H^\infty(G)}^p,
$$
where $G = \varphi(\D)$ is a John domain and $R$ is a rational function of degree at most $n$ with the poles outside
$\overline{G}$.  However, we do not know whether this inequality holds true.  
\medskip
\\
{\bf Proof.} Put $M = \|R\|_{H^\infty(G)}$.
Let us make the change of the variable $\zeta = \varphi(z)$. Taking into account that 
$d_G(\zeta) \le |\varphi'(z)|(1-|z|^2)$, we obtain
$$
I_{p,\beta}(R) \lesssim \int_{\D} |(R\circ\varphi)'(z)|^p |\varphi'(z)|^{2-p+\beta} (1-|z|)^{\beta}\, dA(z).
$$

Statement 3 follows from Theorem 2. Indeed, $1<p-\beta \le 2$ and, using the inequality 
$|(R\circ\varphi)'(z)| (1-|z|) \le M$, we obtain  
$$
I_{p,\beta}(R) \le M^\beta \int_{\D} |(R\circ\varphi)'(z)|^{p-\beta} |\varphi'(z)|^{2-(p-\beta)}\, dA(z) 
\lesssim n^{p-\beta-1}M^p.
$$ 

Let us prove statement 1. The quantity $d_G(z)$ is bounded, so it suffices to prove the statement
for $\beta \in (p-1, p-2 + \gamma]$. 
For such $\beta$ it follows from the inequality $p-\beta <1$ and the inclusion 
$\varphi' \in H^\gamma \subset H^{2-p+\beta}$ that
$$
I_{p,\beta}(R) \lesssim M^p \int_0^1 \bigg(\int_0^{2\pi}|\varphi'(re^{it})|^{2-p+\beta} dt\bigg) \frac{r dr}{(1-r)^{p-\beta}} \lesssim M^p.
$$
At the first step we used the inequality $|(R\circ\varphi)' (z) |(1-|z|) \le M$.

Consider the most interesting statement 2:  $\beta = p-1$. Let $1\le p <2$. Put $s=1-\frac{1}{n}$.
Applying the H\"older inequality with exponents $\frac{2}{p}$ and $\frac{2}{2-p}$ and inequality \eqref{area}, we get
$$
\begin{aligned}
	\int_{\{|z|\le s\}} & |(R\circ\varphi)'(z)|^p |\varphi'(z)| (1-|z|)^{p-1}\, dA(z)   \\
	& \leq \bigg( \int_{\{|z|\le s\}} 
	|(R\circ\varphi)'(z)|^2 (1-|z|)\, dA(z) \bigg)^{\frac{p}{2}}
	\bigg( \int_{\{|z|\le s\}} \frac{|\varphi'(z)|^{\frac{2}{2-p}}}{1-|z|}\, dA(z) \bigg)^{1-\frac{p}{2}} \\
	& \lesssim M^p \bigg( \int_0^s \bigg(\int_0^{2\pi} |\varphi'(re^{it})|^{\frac{2}{2-p}} dt \bigg) \frac{r dr}{1-r} \bigg)^{1-\frac{p}{2}} 
	\lesssim \big( \log n\big)^{1- \frac{p}{2}} M^p.
\end{aligned}
$$
It remains to estimate the integral over the set $\{s<|z| <1\}$:
$$
	\int_{\{s<|z| <1\}}  |(R\circ\varphi)'(z)|^p |\varphi'(z)| (1-|z|)^{p-1}\, dA(z)
	$$
	$$
	  \le 
 M^{p-1} 
	\int_{\{s<|z| <1\}} |(R\circ\varphi)'(z)| \cdot |\varphi'(z)|\, dA(z) 
	$$$$
	 \lesssim M^{p-1} \bigg(\int_{\{s<|z| <1\}} |(R\circ\varphi)'(z)|^2 dA(z)\bigg)^{1/2}  \bigg(\int_{\{s<|z| <1\}}
	|\varphi'(z)|^2 dA(z)\bigg)^{1/2} \lesssim M^p.
$$
In the last inequality we used the fact that 
$$
\int_{\{s<|z| <1\}} |(R\circ\varphi)'(z)|^2 dA(z) \lesssim n M^2,
$$
since $R\circ\varphi$ covers the disk of radius $M$ with multiplicity at most $n$, as well as the inclusion 
$\varphi'\in H^2$.

The case $p\ge 2$ is trivial:
$$
I_{p, p-1}(R) \lesssim M^{p-2} \int_\D |(R\circ\varphi)'(z)|^2 (1-|z|)\, dA(z)\lesssim M^p,
$$
i.e., the quantity $I_{p, p-1}(R)$ is uniformly bounded over $R$ and $n$. 

Theorem 5 is proved.

\end{document}